\documentclass[12pt,reqno]{amsart}
\usepackage{geometry}                
\usepackage{amsmath,amssymb}
\usepackage{graphicx}
\usepackage{amssymb,latexsym, amsmath}
\usepackage{amsfonts,amsthm}
\usepackage{amscd}
\usepackage{mathrsfs}
\usepackage{hyperref}
\usepackage{eucal}
\usepackage{epstopdf}
\usepackage{amsgen}
\usepackage{xspace}
\usepackage{verbatim}
\usepackage{stmaryrd}
\usepackage{enumitem}
\newlist{steps}{enumerate}{1}
\setlist[steps, 1]{label = Step \arabic*:}

\newcommand{\Gd}{\Delta}
\newcommand{\gd}{\delta}

\newcommand{\inpt}[1]{\langle #1 \rangle}

\newcommand{\ga}{\gamma}

\newcommand{\gl}{\lambda}

\newcommand{\pdr}{\partial}

\newcommand{\beq}{\begin{equation}}
\newcommand{\eeq}{\end{equation}}
\newcommand{\bea}{\begin{align}}
\newcommand{\eea}{\end{align}}
\newcommand{\bthm}{\begin{theorem}}
\newcommand{\ethm}{\end{theorem}}
\newcommand{\bpr}{\begin{proof}}
\newcommand{\epr}{\end{proof}}
\newcommand{\bcl}{\begin{corollary}}
\newcommand{\ecl}{\end{corollary}}
\newcommand{\bpn}{\begin{proposition}}
\newcommand{\epn}{\end{proposition}}
\newcommand{\bre}{\begin{remark}}
\newcommand{\ere}{\end{remark}}
\newcommand{\bdf}{\begin{definition}}
\newcommand{\edf}{\end{definition}}
\newcommand{\bss}{\begin{align*}}
\newcommand{\ess}{\end{align*}}

\newcommand{\bl}{\label}

\newcommand{\mR}{\mathbb{R}}

\newcommand{\mZ}{\mathbb{Z}}

\newtheorem{theorem}{Theorem}[section]
\newtheorem{corollary}[theorem]{Corollary}

\newtheorem{proposition}[theorem]{Proposition}

\theoremstyle{definition}
\newtheorem{definition}[theorem]{Definition}
\theoremstyle{remark}
\newtheorem{remark}{Remark}

\numberwithin{equation}{section}

\begin{document}

\title[Synchronization of Cellular Neural Networks]{Feedback Synchronization of FHN Cellular Neural Networks}

\author[L. Skrzypek]{Leslaw Skrzypek}
\address{Department of Mathematics and Statistics, University of South Florida, Tampa, FL 33620, USA}
\email{skrzypek@usf.edu}
\thanks{}

\author[Y. You]{Yuncheng You}
\address{Department of Mathematics and Statistics, University of South Florida, Tampa, FL 33620, USA}
\email{you@mail.usf.edu}
\thanks{}


\date{June 18, 2020}



\begin{abstract} 
In this work we study the synchronization of ring-structured cellular neural networks modeled by the lattice FitzHugh-Nagumo equations with boundary feedback. Through the uniform estimates of solutions and the analysis of dissipative dynamics, the synchronization of this type neural networks is proved under the condition that the boundary gap signal exceeds the adjustable threshold. 
\end{abstract}

\maketitle

\section{\textbf{Introduction}}

Cellular neural network (briefly CNN) was invented by Chua and Yang \cite{Chua1, Chua2}. CNN consists of network-like interacted processors of continuous-time analog or digital signals. In some sense the dynamics of CNN can be analyzed as a lattice dynamical system oftentimes generated by lattice differential equations in time \cite{Chua3, CR, S2}. 

From the theoretical viewpoint, the theory of cellular neural networks and, more recently, of convolutional neural networks (also called CNN) and variants of complex neural networks is closely linked to discrete nonlinear partial differential equations and delay differential equations as well as the spatially discrete Fourier transform and constrained optimization \cite{Chow1, Chow2, Chua3}. In the dramatically broadened application fronts to machine leaning, deep learning, and general artificial intelligence, the most prominent area is the image processing especially in medical visualization techniques \cite{Chua2, CR, S1}. 

Many complicated computational problems can be formulated as multi-layer and parallel tasks for processing signal values on a geometric grid with direct interaction and transmission in a local neighborhood. The cloning template of each cell model and the coupling design mimic biological pattern formation in the brain and nerves are the two features of CNN. The concepts and models of CNN are based on some aspects of neurobiology \cite{S2} in terms of sheet-like layers or arrays of massively interconnected excitable neurons and the implementation counterpart by electronic integrated circuit chips on the other hand.

Rapidly expanding applications of CNN and all kinds of complex neural networks \cite{CR, Hp, I, W} stimulate frontier and interdisciplinary researches from data analysis to mathematical and statistical modeling to deep learning algorithms, softwares, and robotics. 

In this work, we consider a CNN model as the discrete replacement of the biological neural networks described by the partial diffusive FitzHugh-Nagumo equations with different coupling patterns \cite{FH, LY, S1, WLZ, Yong},
\begin{align*}
	\frac{\pdr u_i}{\pdr t} &= a \Gd u_i + f(u_i) - b w_i, \\
	\frac{\pdr w_i}{\pdr t} &= c u_i - \gd w_i. 
\end{align*}
We shall use the discrete Laplacian templates to approximate the Laplacian operator on a 1D domai \cite{Chow1, S2}. For the prototype FitzHugh-Nagumo equations, the nonlinear term is
$$
	f(s) = s(s-\alpha)(1 - s)
$$ 
and the parameters $0 < \alpha < 1$ and $0 < c \ll 1$. 

Consider a layer of cellular neural network, which consists of 1D cells at the grid points $\{ih: i = 1, 2, \cdots , n\}$ for $h > 0$, in a ring-structure. We shall study synchronization dynamics of the FitzHugh-Nagumo lattice equations: 
\beq \bl{CN}	
	\begin{split}
	\frac{dx_i}{dt} &= a (x_{i-1} - 2x_i + x_{i+1}) + f(x_i) - b y_i + p u_i, \quad 1 \leq i \leq n, \\
	\frac{dy_i}{dt} &= c x_i - \gd y_i, \quad 1 \leq i \leq n,
	\end{split}
\eeq
where $t > 0$, the integer $n \geq 4$, and the discrete Laplacian operator
$$
	D_i (x) = a (x_{i+1} - 2x_i + x_{i-1}),
$$
can be called the synaptic law of cell coupling. In this model of CNN, we impose the periodic boundary condition
\beq \bl{pbc}
	x_0(t) = x_n(t), \quad x_{n+1}(t) = x_1(t)
\eeq
and the boundary feedback control $\{u_i\}_{i=1}^n$, 
\beq \bl{bfc}
	\begin{split}
	&u_1(t) = u_{n+1}(t) = x_n (t) - x_1(t),     \\
	&u_i (t) = 0, \quad 2 \leq i \leq n - 1, \\
	&u_{n}(t) = u_0(t) = x_1(t) - x_n(t).
	\end{split}
\eeq
for $t \geq 0$, where $p > 0$ is the controllable feedback constant and $x_n(t) - x_1(t)$ measures the boundary gap signal between the two endpoints of the cellular neural network. The initial conditions for the system \eqref{CN} are denoted by
\beq \bl{inc}
	x_i(0) = x_i^0 \in \mR \quad \text{and} \quad y_i (0) = y_i^0 \in \mR, \qquad 1 \leq i \leq n.
\eeq
All the parameters in this system \eqref{CN} are positive constants. 

We make the following Assumption: The scalar function $f \in C^1 (\mR, \mR)$ satisfies 
\beq \bl{Asp}
	\begin{split}
	&f(s) s \leq  - \gl s^4 + \beta, \quad s \in \mathbb{R}, \\
	& f^{\,\prime} (s) \leq \ga, \quad s \in \mathbb{R}, \\
	\end{split}
\eeq
where $\gl, \beta$ and $\ga$ are positive constants. The typical nonlinearity $f(s) = s(s - \alpha)(1 - s)$ shown above in the FitzHugh-Nagumo model satisfies the Assumption \eqref{Asp}:
\begin{align*}
	f(s)s &= - \alpha s^2 + (\alpha + 1)s^3- s^4 \leq -\alpha s^2 + \left(\frac{1}{2}s^4 + 2^3 (\alpha +1)^4\right) - s^4 \\
	&\leq - \left(\alpha s^2 + \frac{1}{2}s^4 \right) + 8(\alpha + 1)^4 \leq - \frac{1}{2} s^4 + 8(\alpha + 1)^4,  \\[3pt]
        f^{\, \prime} (s) &= - \alpha + 2(\alpha +1)s - 3s^2 \leq - \alpha + (\alpha +1)^2 - 2s^2 \leq 1 + \alpha + \alpha^2.
\end{align*}


Synchronization plays a significant role for biological neural networks and for the artificial neural networks as well. Fast synchronization may lead to enhanced functionality and performance of complex neural networks. 

In recent years, the dynamical behavior and problems of complex and large-scale networks including convolutional neural networks in machine learning and deep learning, Internet networks, epidemic spreading networks, and social networks attract many interdisciplinary research interests, cf.\cite{A, I, W}. Synchronization for the CNN modeled by PDE, delay differential equations, or lattice differential equations is one of the essential topics in the theoretical analysis of artificial intelligence. 

Synchronization for biological neural networks has been studied by several mathematical models and methods. This topic has been studied for the diffusive FitzHuigh-Nagumo networks of neurons coupled by clamped gap junctions \cite{AA, A, AC, IJ, Yong}, the mean field couplings of Hodgkin-Huxley and FitzHuigh-Nagumo neuron networks  \cite{Dick, QT}, and the chaotic neural networks and stochastic neural networks \cite{Dick, SG}.

Recently we proved results on the exponential synchronization of the boundary coupled Hindmarsh-Rose neural networks in \cite{PLY, PY} and the boundary coupled partly diffusive FitzHugh-Nagumo neural networks in \cite{LY}.

The feature of this work is to provide a sufficient condition for realization of the feedback synchronization of the proposed FitzHugh-Nagumo (FHN) cellular neural networks with boundary control. The quantitative threshold condition for synchronization is explicitly expressed in terms of the parameters and can be adjusted by the feedback strength coefficient $p$ in  applications. 

\section{\textbf{Uniform Estimates and Dissipative Dynamics}}

Define the following Hilbert space:
$$
	H = \ell^2 (\mZ_n, \mR^{2n}) = \{ (x, y) = \{(x_i, y_i): 1 \leq i \leq n \}\} 
$$
where $\mZ_n = \{1 ,2, \cdots, n\}$ and $n \geq 4$. The norm in $H$ is denoted and define by $\|(x, y)\|^2 = \|x\|^2 + \|y\|^2 =  \sum_{i=1}^n (| x_i  |^2 + |y_i |^2)$.
The inner-product of $H$ or $\mathbb{R}^n$ is denoted by $\inpt{\,\cdot , \cdot\,}$.

Since there exists a unique local solution in time of the initial value problem \eqref{CN}-\eqref{inc} under the Assumption \eqref{Asp} that the right-side functions in \eqref{CN} are locally Lipschitz continuous, in this section we shall first prove the global existence in time of the solutions in the space $H$. By the uniform estimates we show the dissipative dynamics of the solution semiflow.

\begin{theorem} \label{Tm}
Under the boundary feedback control \eqref{bfc}, for any given initial state $(x^0, y^0) = ((x_1^0, y_1^0), \cdots , (x_n^0, y_n^0)) \in H$, there exists a unique solution, 
$(x(t), y(t)) = ((x_1 (t, x_1^0), y_1(t, y_1^0)), \cdots, (x_n (t, x_n^0)), y_n(t, y_n^0)), t \in [0, \infty)$,  of the initial value problem \eqref{CN}-\eqref{inc} for this cellular neural network. 
\end{theorem}

\begin{proof}
Multiply the $x_i$-equation in \eqref{CN} by $C_1 x_i(t)$ for $1 \leq i \leq n$, where the constant $C_1 > 0$ is to be chosen, then sum them up and by the Assumption \eqref{Asp} to get	
\begin{equation} \bl{u1}
	\begin{split}
	&\frac{C_1}{2} \frac{d}{dt} \sum_{i = 1}^n |x_i |^2 = C_1 \sum_{i=1}^n \left[ a(x_{i-1} -2x_i + x_{i+1}) x_i + f(x_i) x_i - b x_i y_i + p u_i x_i \right]   \\
	\leq &\, C_1 \sum_{i=1}^n a(x_{i-1} - 2x_i+ x_{i+1}) x_i  \\
	+ &\, C_1 \sum_{i=1}^n \left[- \gl |x_i |^4 + \beta + \frac{b}{2}\, |x_i|^2 + \frac{b}{2}\, |y_i |^2 \right] - C_1 p(x_1 - x_n)^2, \quad t \in I_{max}, 
	\end{split}
\end{equation}
where $I_{max} = [0, T_{max})$ is the maximal existence interval of the solution. According to the discrete "divergence" formula and $x_0(t) = x_n(t), \, x_{n+1}(t) = x_1(t)$ due to the periodic boundary condition \eqref{pbc}, we have
\beq \bl{key}
	\begin{split}
	&\sum_{i=1}^n (x_{i-1} -2x_i + x_{i+1})x_i = \sum_{i=1}^n (x_{i+1} - x_i)x_i - \sum_{i=1}^n (x_i - x_{i-1})x_i   \\
	= &\, \left[\sum_{i=1}^{n-1} (x_{i+1} - x_i)x_i - \sum_{i=2}^n (x_i - x_{i-1})x_i\right] + (x_{n+1} - x_n)x_n - (x_1 - x_0)x_1 \\
	= &\, - \sum_{i=2}^n (x_i - x_{i-1})^2 - (x_1 - x_0)^2 = - \sum_{i=1}^{n} (x_{i} - x_{i-1})^2 \leq 0.
	\end{split}
\eeq
Then \eqref{u1} with \eqref{key} yields the differential inequality
\beq \bl{u2}
	\begin{split}
	&C_1 \frac{d}{dt} \sum_{i = 1}^n |x_i (t)|^2 + 2C_1 \left[\sum_{i=2}^{n} a(x_{i} - x_{i-1})^2 + p (x_1 - x_n)^2\right]  \\[2pt]
	\leq &\,C_1 \sum_{i=1}^N \left[- 2\gl |x_i (t)|^4 + 2\beta + b |x_i(t)|^2 + b |y_i (t)|^2 \right],  \quad t \in I_{max}.
	\end{split}
\eeq

Next multiply the $y_i$-equation in \eqref{CN} by $y_i(t)$ for $1 \leq i \leq n$ and then sum them up. By using Young's inequality, we obtain
\begin{equation} \bl{w1}
	\begin{split}
	\frac{1}{2} \frac{d}{dt} \sum_{i=1}^n |y_i(t) |^2 &= \sum_{i=1}^n ( cx_i y_i - \gd y_i^2) \leq \sum_{i=1}^n \left[\left(\frac{c^2}{\gd} x_i^2 + \frac{1}{4} \gd \,y_i^2\right) - \gd \,y_i^2\right]  \\[3pt]
	= &\, \sum_{i=1}^n \left[\frac{c^2}{\gd} \, |x_i (t)|^2 - \frac{3}{4} \gd \,|y_i(t)|^2\right], \quad \text{for} \;\, t \in I_{max}.
	\end{split}
\end{equation}
Now add the above two inequalities \eqref{u2} and doubled \eqref{w1}. We obtain
\beq \bl{uw}
	\begin{split}
        &\frac{d}{dt} \sum_{i = 1}^n \left(C_1 |x_i (t)|^2 + | y_i (t)|^2 \right)  + 2C_1 \left[\sum_{i=1}^{n} a(x_{i} - x_{i-1})^2 + p (x_1 - x_n)^2\right]  \\[2pt]
        \leq &\sum_{i=1}^n \left[ \left(C_1 b + \frac{2c^2}{\gd}\right) |x_i(t)|^2 - 2C_1 \gl |x_i (t)|^4 + 2C_1 \beta \right] \\[2pt]
        + &\sum_{i=1}^n \left[\left(C_1 b - \frac{3 \gd}{2}\right) |y_i (t)|^2 \right], \quad t \in I_{mac} = [0, T_{max}).
	\end{split}
\eeq
We choose constant 
\beq \bl{C1}
	C_1 = \frac{\gd}{2b} \quad \text{so that} \quad C_1 b - \frac{3\gd}{2} = - \gd.
\eeq
Then from \eqref{uw} with the fact $2C_1 [ \cdots ] \geq 0$ on the left-hand side and from the choice \eqref{C1}, we have
$$
	 \frac{d}{dt} \sum_{i = 1}^n \left(C_1 |x_i |^2 + | y_i |^2 \right) \leq \sum_{i=1}^n \left[ \left(C_1 b + \frac{2c^2}{\gd}\right)|x_i(t)|^2- 2C_1 (\gl |x_i (t)|^4 + \beta) - \gd |y_i (t)|^2\right] 
$$
and consequently,

\beq \bl{Cuw}
	\begin{split}
	&\frac{d}{dt} \sum_{i = 1}^n \left(C_1 |x_i(t) |^2 + | y_i (t)|^2 \right) + \gd  \sum_{i = 1}^n \left( C_1 |x_i (t) |^2 + | y_i (t)|^2 \right)   \\
	\leq &\, \sum_{i = 1}^n \left[ \left(C_1 b + C_1 \gd + \frac{2c^2}{\gd}\right)|x_i(t)|^2- 2C_1 (\gl |x_i (t)|^4 + \beta)\right]   \\
	= &\, \sum_{i = 1}^n  \left[ \left(\frac{\gd^2}{2b} + \frac{\gd}{2} + \frac{2c^2}{\gd}\right) |x_i(t)|^2- \frac{\gd \gl}{b} |x_i (t)|^4 + \frac{\gd \beta}{b} \right], \quad t \in I_{max}.
	\end{split}
\eeq
Completing square shows that
\begin{align*}
	 &\left(\frac{\gd^2}{2b} + \frac{\gd}{2} + \frac{2c^2}{\gd} \right) |x_i(t)|^2- \frac{\gd \gl}{b} |x_i (t)|^4  \\
	 = &\, - \frac{\gd \gl}{b} \left[ | x_i (t)|^2 - \frac{b}{2\gd \gl} \left(\frac{\gd^2}{2b} + \frac{\gd}{2} + \frac{2 c^2}{\gd}\right) \right]^2 + C_2
\end{align*}
and
\beq \bl{C2}
	C_2 = \frac{b}{4\gd \gl} \left(\frac{\gd^2}{2b} + \frac{\gd}{2} + \frac{2 c^2}{\gd}\right)^2.
\eeq
Therefore, \eqref{Cuw} yields
\beq \bl{Suw}
	\frac{d}{dt} \sum_{i = 1}^n \left(C_1 |x_i |^2 + | y_i |^2 \right) + \gd  \sum_{i = 1}^n \left( C_1 |x_i |^2 + | y_i |^2 \right) \leq n \left(C_2 + \frac{\gd \beta}{b}\right), \quad t \in I_{max}.
\eeq
Apply the Gronwall inequality to \eqref{Suw}. Then we get the following bounded estimate for all the solutions of the system of equations \eqref{CN}-\eqref{inc},
\beq \label{dse}
	\begin{split}
	& \sum_{i=1}^n \left(|x_i (t, x_i^0) |^2 + |y_i (t, y_i^0)|^2 \right) \\
	\leq &\, \frac{1}{\min \{C_1, 1\}} \left[e^{- \gd \, t} \sum_{i=1}^n \left(C_1 |x_i^0 |^2 + |y_i^0 |^2 \right) + \frac{n}{\gd} \left(C_2 + \frac{\gd \beta}{b}\right)\right], \;\; t \in [0, \infty). \\
	\end{split}
\eeq
Here it is shown that $I_{max} = [0, \infty)$ for all the solutions because they will never blow up at any finite time.  Thus it is proved that for any given initial state there exists a unique global solution $((x_1 (t, x_1^0), y_1(t, y_1^0)), \cdots, (x_n (t, x_n^0)), y_n(t, y_n^0))$ in $H$. 
\end{proof}

The global existence and uniqueness of the solutions to the initial value problem \eqref{CN}-\eqref{inc} and their continuous dependence on the initial data enable us to define the solution semiflow $\{S(t): H \to H\}_{t \geq 0}$ of this system of the FitzHugh-Nagumo cellular neural network:
$$
	S(t): ((x_1^0, y_1^0) \cdots, (x_n^0, y_n^0)) \longmapsto ((x_1(t, x_1^0), y_1 (t, y_1^0)), \cdots, (x_n(t, x_n^0), y_n (t, y_n^0))).
$$
We call $\{S(t)\}_{t \geq 0}$ the semiflow of the FitzHugh-Nagumo CNN.

\begin{theorem} \label{Dsp}
	The semiflow $\{S(t)\}_{t \geq 0}$ of the FitzHugh-Nagumo CNN in the space $H$ is dissipative in the sense that there exists a bounded ball 
\beq \label{abs}
	B^* = \{g \in H: \| g \|^2 \leq Q\}
\eeq 
where the constant 
\beq \bl{Q}
	Q = \frac{1}{\min \{C_1, 1\}} \left[1 +  \frac{n}{\gd} \left(C_2 + \frac{\gd \beta}{b}\right) \right]
\eeq
such that for any given bounded set $B \subset H$, there is a finite time $T_B > 0$ and all the solutions with the initial state inside the set $B$ will permanently enter the ball $B^*$ for $t \geq T_B$.
\end{theorem}

\begin{proof}
The uniform estimate \eqref{dse} implies that
	\beq \label{lsp}
	\limsup_{t \to \infty} \,\sum_{i=1}^n \left(|x_i (t, x_i^0) |^2 + |y_i (t, y_i^0)|^2 \right) < Q 
	\eeq
for all solutions of \eqref{CN} with any initial data $((x_1^0, y_i^0), \cdots, (x_n^0, y_n^0)) \in H$. Indeed for any given bounded set $B = \{g \in H: \|g \|^2 \leq \rho \}$ in $H$, there is a finite time 
$$
	T_B = \frac{1}{\gd} \log^+ (\rho \, \max \{C_1, 1\} )
$$	
such that 
$$
	e^{- \gd \, t} \sum_{i=1}^n \left(C_1 |x_i^0 |^2 + |y_i^0 |^2 \right) < 1, \quad \text{for} \;\, t \geq T_B,
$$
which means all the solution trajectories started from the set $B$ will permanently enter the bounded ball $B^*$ shown in \eqref{abs} for $t \geq T_B$. Therefore, this semiflow is dissipative.
\end{proof}

\section{\textbf{Synchronization of the FitzHugh-Nagumo CNN}} 

Define the differences of solutions for two adjacent indexed cells of the FitzHugh-Nagumo CNN \eqref{CN} to be
$$
	V_i (t) = x_i(t) - x_{i-1}(t), \quad W_i (t) = y_i(t) - y_{i-1}(t), \quad \text{for}\;\,  i = 1, \cdots, n. 
$$
Consider the system of the \emph{differencing} equations for this CNN. For $i = 1, \cdots, n$,
\beq \bl{dHR}
	\begin{split}
		\frac{\pdr V_{i}}{\pdr t} & = a (V_{i-1} - 2V_{i} + V_{i+1})+ f(x_i) - f(x_{i-1}) - b W_{i} + p(u_i - u_{i-1}),  \\
		\frac{\pdr W_{i}}{\pdr t} & = c\, V_{i} - \gd W_{i}.
	\end{split}
\eeq
The periodic boundary condition $V_0 (t) = V_n (t), \, V_{n+1}(t) = V_1 (t)$ holds due to \eqref{pbc}. 

Here is the main result on the feedback synchronization of the proposed FitzHugh-Nagumo cellular neural networks.

\begin{theorem} \bl{ThM}
	If the following threshold condition for the boundary gap signal of the FitzHugh-Nagumo cellular neural network \eqref{CN}-\eqref{bfc} is satisfied,  
\beq \bl{SC}
	 \liminf_{t \to \infty} \, (x_n(t) - x_1(t))^2 > \left(1 + \frac{1}{p} \left(\gd + \ga + |c - b|\right)\right)Q, 
\eeq
where the constant $Q > 0$ is given in \eqref{Q}, then this cellular neural network is asymptotically synchronized in the space $H$ at a uniform exponential rate.
\end{theorem}

\begin{proof}
	 Multiply the first equation in \eqref{dHR} by $V_{i} (t)$ and the second equation in \eqref{dHR} by $W_{i} (t)$. Then sum them up for all $1 \leq i \leq n$ and use the Assumption  \eqref{Asp} to get
\beq \bl{eG}
	\begin{split}
	&\frac{1}{2} \frac{d}{dt} \sum_{i=1}^n \left(|V_{i}|^2 + |W_{i}|^2\right) - \sum_{i=1}^n a (V_{i-1} - 2V_i + V_{i+1})V_i \\
	= &\, \sum_{i=1}^n \left[ (f(x_i) - f(x_{i-1})) V_{i} + (c - b)V_{i} W_{i} - \gd \,|W_{i}|^2 + p(u_i - u_{i-1})V_i \right] \\
	\leq &\,\sum_{i=1}^n \left[f^{\,\prime} (\xi x_i + (1-\xi) x_{i-1} ) V^2_{i} + (c - b)V_{i} W_{i} - \gd \,|W_{i}|^2 + p(u_i - u_{i-1})V_i \right]   \\
	\leq &\, \sum_{i=1}^n \left[\ga |V_{i}|^2 + |c - b |(|V_{i}|^2 + |W_{i}|^2) - \gd \,|W_{i}|^2 + p(u_i - u_{i-1})V_i \right],
	\end{split}
\eeq
where $ 0 \leq \xi \leq 1$. By the periodic boundary condition for the differencing equations \eqref{dHR} we have 
\begin{align*}
	&- \sum_{i=1}^n a (V_{i-1} - 2V_i + V_{i+1})V_i = - \sum_{i=1}^n a (V_{i+1} - V_i)V_i + \sum_{i=1}^n a (V_i - V_{i-1})V_i   \\
	= &\,- \left[\sum_{i=1}^{n-1} a(V_{i+1} - V_i)V_i - \sum_{i=2}^n a(V_i - V_{i-1})V_i\right] - a(V_{n+1} - V_n)V_n + a(V_1 - V_0)V_1 \\
	= &\, \sum_{i=2}^{n} a(V_i - V_{i-1})^2 + a(V_1^2 + V_n^2) - a(V_{n+1}V_n + V_0 V_1)   \\
	= &\, \sum_{i=2}^{n} a(V_i - V_{i-1})^2 + a(V_1^2 + V_0^2) - 2a V_1 V_0 \\
	= &\, \sum_{i=2}^{n} a(V_i - V_{i-1})^2 + a(V_1^2 - V_0^2) = \sum_{i=1}^{n} a(V_i - V_{i-1})^2 \geq 0.
\end{align*}
From the above two inequalities, we obtain
\beq \bl{VW}
	\begin{split}
	&\frac{1}{2} \frac{d}{dt} \sum_{i=1}^n \left(|V_i (t)|^2 + |W_i (t)|^2\right)   \\
	\leq &\, \sum_{i=1}^n \left[\ga |V_{i}|^2 + |c - b |(|V_{i}|^2 + |W_{i}|^2) - \gd \,|W_{i}|^2 + p(u_i - u_{i-1})V_i \right].
	\end{split}
\eeq
The boundary feedback \eqref{bfc} and the perioc boundary condition \eqref{pbc} infer that
\beq \bl{pu}
	\begin{split}
	 & \sum_{i=1}^n \, p(u_i - u_{i-1})V_i = \sum_{i=1}^n \, p(u_i - u_{i-1})(x_i - x_{i-1})  \\
	 = &\, p \left[(u_1 - u_0)(x_1 - x_0) + (u_2 - u_1)(x_2 - x_1) + (u_n - u_{n-1})(x_n - x_{n-1})\right] \\[4pt]
	 = &\, p \left[(u_1 - u_n)(x_1 - x_n) - u_1(x_2 - x_1) + u_n(x_n - x_{n-1})\right]  \\[4pt]
	 = &\, p \left[2(x_n - x_1)(x_1 - x_n) - (x_n - x_1)(x_2 - x_1) + (x_1 - x_n)(x_n - x_{n-1})\right]  \\[2pt]
	 = &\, p \left[- 2(x_n - x_1)^2 + (x_n - x_1)(x_1 - x_2 + x_{n-1} - x_n)\right]  \\[2pt]
	 = &\, p \left[- 3(x_n - x_1)^2 + (x_n - x_1)(x_{n-1} - x_2)\right]  \\[2pt]
	 \leq &\, p \left[- 2(x_n - x_1)^2 + (x_{n-1} - x_2)^2 \right].
	 \end{split}
\eeq
Substitute \eqref{pu} into \eqref{VW}. Then we get the following differential inequality
\begin{equation*}
	\begin{split}
	&\frac{d}{dt} \sum_{i=1}^n \left(|V_i (t)|^2 + |W_i (t)|^2\right) + 4p (x_n(t) - x_1(t))^2 \\
	\leq &\,\sum_{i=1}^n 2\left[\ga |V_{i}|^2 + |c - b |(|V_{i}|^2 + |W_{i}|^2) - \gd \,|W_{i}|^2\right] + 2p(x_{n-1}(t) - x_2(t))^2.
	\end{split}
\end{equation*}
Hence it holds that
\begin{equation}  \bl{Kq}
	\begin{split}
	&\frac{d}{dt} \sum_{i=1}^n \left(|V_i (t)|^2 + |W_i (t)|^2\right) + 2\gd \sum_{i=1}^n (|V_i(t)|^2 + |W_i(t)|^2) + 4p (x_n(t) - x_1(t))^2 \\
	\leq &\,\sum_{i=1}^n 2\left[(\gd + \ga)|V_i(t)|^2 + |c - b |(|V_{i}|^2 + |W_{i}|^2) \right] + 2p(x_{n-1}(t) - x_2(t))^2, 
	\end{split}
\end{equation}
for $t > 0$. Note that \eqref{lsp} in Theorem \ref{Dsp} confirms that for all solutions of \eqref{CN},
$$
	\limsup_{t \to \infty}\, \sum_{i=1}^n \left(|x_i(t, x_i^0)|^2 + |y_i(t, y_i^0|^2\right) < Q.
$$
Thus for any given bounded set $B \subset H$ and any initial data $((x_1^0, y_i^0), \cdots , (x_n^0, y_n^0)) \in B$, there is a finite time $T_B \geq 0$ such that 
\beq \bl{bd}
	\begin{split}
	&\sum_{i=1}^n 2\left[(\gd + \ga)|V_i(t)|^2 + |c - b |(|V_{i}|^2 + |W_{i}|^2)\right] + 2p(x_{n-1}(t) - x_2(t))^2   \\
	< &\, 4 \left(\gd + \ga + |c - b| \right) Q + 4p\,Q = 4 \left(\gd + \ga + |c - b|  + p\right) Q, \quad \text{for} \;\,   t \geq T_B.
	\end{split}
\eeq
Combining \eqref{Kq} and \eqref{bd}, we have shown that 
\beq \bl{Mq}
	\begin{split}
	\frac{d}{dt} \sum_{i=1}^n \left(|V_i (t)|^2 + |W_i (t)|^2\right) &+ 2\gd \sum_{i=1}^n (|V_i(t)|^2 + |W_i(t)|^2) + 4p (x_n(t) - x_1(t))^2   \\
	&< 4 \left(\gd + \ga + |c - b|  + p\right) Q, \quad \text{for}  \;\;  t \geq T_B.
	\end{split}
\eeq

	Under the threshold condition \eqref{SC} of this theorem, for any given initial state $(x^0, y^0) = ((x_1^0, y_1^0), \cdots, (x_n^0, y_n^0)) \in H$ as a set $B$ of single point, there exists a finite time $T_{(x^0, \, y^0)}  > 0$ such that the differential inequality \eqref{Mq} holds for $t > T_{(x^0, \, y^0)}$ and 
$$
	\sum_{i=1}^n \left(|x_i (t, x_i^0) |^2 + |y_i (t, y_i^0)|^2 \right) < Q, \quad \text{for} \;\; t \geq T_{(x^0, \, y^0)}.
$$
Moreover,
$$
	(x_n(t) - x_1(t))^2 > \left(1 + \frac{1}{p} \left(\gd + \ga + |c - b|\right)\right)Q,  \quad \text{for} \;\; t \geq T_{(x^0, \, y^0)},
$$
so that
\beq \bl{pl}
	p (x_n(t) - x_1(t))^2 > \left(\gd + \ga + |c - b| + p \right)Q,  \quad \text{for} \;\; t \geq T_{(x^0, \, y^0)}.
\eeq
It follows from \eqref{Mq} and \eqref{pl} that
\beq \bl{Gwq}
	\frac{d}{dt} \sum_{i=1}^n \left(|V_i (t)|^2 + |W_i (t)|^2\right) + 2\gd \sum_{i=1}^n (|V_i(t)|^2 + |W_i(t)|^2) < 0, \;\; \text{for} \;\, t \geq T_{(x^0, \, y^0)}.
\eeq
Finally, the Gronwall inequality applied to \eqref{Gwq} shows that
\beq \bl{Syn}
	\begin{split}
	\sum_{i=1}^n (|V_i(t)|^2 + |W_i(t)|^2) &\leq e^{- \gd (t - T_{(x^0, \, y^0)})} \sum_{i=1}^n (|V_i (T_{(x^0, \, y^0)})|^2 + |W_i (T_{(x^0, \, y^0)})|^2)  \\[2pt]
	&\leq 2 e^{- \gd (t - T_{(x^0, \, y^0)})}\,Q \to 0,  \quad \text{as} \;\; t \to \infty.
	\end{split}
\eeq
Then it is proved that for all solutions of the problem \eqref{CN}-\eqref{bfc} for this FitzHugh-Nagumo CNN with the boundary feedback, 
\beq \bl{gij}
	\lim_{t \to \infty} \sum_{i = 1}^n \left(|(x_i (t, x_i^0) - x_{i-1}(t, x_{i-1}^0))|^2 + |(y_i (t, y_i^0) - y_{i-1}(t, y_{i-1}^0))|^2 \right) = 0.
\eeq
This FHN cellular neural network with boundary feedback is asymptotically synchronized in the space $H$ at a uniform exponential rate. 
\end{proof}

This result provides a sufficient condition for feedback synchronization of the FitzHugh-Nagumo complex neural networks with boundary control. The threshold condition \eqref{SC} needs to be satisfied by the boundary gap signal $\liminf_{t \to \infty}\, (x_n(t) - x_1(t))^2$ between the two boundary cells. And the threshold in \eqref{SC} is adjustable by the designed feedback coefficient $p$ in applications. 

\bibliographystyle{amsplain}

\end{document}